\documentclass[conference]{IEEEtran}
\IEEEoverridecommandlockouts
\usepackage{amsmath,amssymb,amsfonts}
\usepackage{graphicx}
\usepackage{tabularx}
\usepackage{pifont}
\usepackage{cite}
\usepackage{bm}
\usepackage{optidef}
\usepackage{flushend}
\usepackage{multirow}
\usepackage{xcolor}
\usepackage[colorlinks=true, linkcolor=blue, citecolor=blue, urlcolor=blue]{hyperref}
\usepackage{url}
\usepackage[figure]{hypcap}

\begin{document}

\title{Incentive-Aligned Vehicle-to-Vehicle Energy Trading via Nash-Integrated Multi-Agent Reinforcement Learning}

\author{\IEEEauthorblockN{Yujin Lin\textsuperscript{1}, Yue Yang\textsuperscript{1}, Hao Wang\textsuperscript{1,2*}}
\IEEEauthorblockA{\textsuperscript{1}Department of Data Science and AI, Faculty of IT, Monash University, Australia \\
\textsuperscript{2}Monash Energy Institute, Monash University, Australia
%\textsuperscript{}Emails: ylin0115@student.monash.edu, yue.yang1@monash.edu, hao.wang2@monash.edu\\
}
\thanks{*Corresponding author: Hao Wang (hao.wang2@monash.edu).}
\thanks{This work was supported in part by the Australian Research Council (ARC) Discovery Early Career Researcher Award (DECRA) under Grant DE230100046.}
}

\maketitle

\begin{abstract}
Vehicle-to-vehicle (V2V) energy trading enables decentralized peer-to-peer energy exchange among electric vehicles (EVs), reducing grid dependency while monetizing surplus capacity. However, coordinating self-interested EV agents with diverse charging needs and uncertain arrival-departure schedules remains challenging. Existing approaches either require centralized optimization with computational limitations or lack fairness guarantees. This paper integrates Nash Bargaining Solution into Multi-Agent Deep Deterministic Policy Gradient, namely Nash-MADDPG, for incentive-aligned V2V energy trading. Nash bargaining determines efficient bilateral pricing, while Nash-guided price proximity rewards align agent learning toward bargaining-optimal strategies. Evaluation over 30-day continuous operation demonstrates an improvement of 61.6\% in social welfare and 62.9\% improvement in trading volume over Double Auction, while achieving superior fairness, such as 40.1\% improvement in Jain's index. Testing across 6-100 agents over a 30-day horizon with continuous vehicle turnover confirms scalability across population size and empirically stable pricing near the Nash Bargaining benchmark.
\end{abstract}


\section{Introduction}
Electric vehicle (EV) adoption challenges power grids during peak charging, causing congestion and voltage instability. Traditional approaches such as centralized charging coordination \cite{wang2023fast}, time-of-use pricing, and vehicle-to-grid (V2G) services rely on centralized infrastructure and focus on EV-grid interactions, leaving direct EV-to-EV energy exchange underutilized.

Vehicle-to-vehicle (V2V) energy trading offers a decentralized alternative where EVs with surplus capacity sell energy directly to EVs needing charge, reducing grid dependency \cite{xu2024vehicle}. However, realizing efficient V2V markets requires solving three intertwined challenges: (i) fair bilateral pricing between self-interested agents with private valuations, (ii) adaptation to dynamic participation as vehicles continuously arrive and depart, and (iii) efficient, individually rational outcomes that encourage voluntary participation.

Existing V2V approaches partially address these challenges. In \cite{khele2023fairness}, V2V interaction was formulated as a centralized fairness-aware optimization, jointly allocating energy and prices for equitable outcomes, but computational cost grows rapidly with agent count. In \cite{kabir2020routing}, routing and scheduling of mobile EV chargers was optimized, though without addressing strategic bidding among stationary EVs. Blockchain-based frameworks address trust: \cite{wang2023fast} developed a blockchain-consensus V2V protocol with Stackelberg pricing where sellers set prices and buyers respond optimally; \cite{kumari2023v2v} proposed wireless V2V sharing with blockchain verification against double-spending. However, these game-theoretic pricing models assume complete information or static conditions that break down when vehicles arrive unpredictably with private valuations. A systematic review of V2V studies \cite{xu2024vehicle} identifies four open challenges (dynamic matching, real-time pricing, incentive alignment, scalability) and concludes that hybrid learning-mechanism approaches are needed.

Nash Bargaining Solution (NBS) \cite{nash1950bargaining} provides axiomatic Pareto efficiency and individual rationality guarantees attractive for energy markets. NBS was first applied to energy trading in \cite{wang2016incentivizing}, pricing bilateral microgrid exchanges at the bargaining point so both parties benefit relative to standalone operation. Subsequent work applied NBS to regional energy systems under wholesale price uncertainty \cite{wang2023nash_energy}, multi-carrier peer-to-peer trading with stochastic programming under renewable uncertainty \cite{alizadeh2024nash_p2p}, and asymmetric settings where agents receive bargaining power proportional to their contribution to shared energy storage \cite{chen2023asymmetric_nash}. However, these approaches are model-based, requiring complete knowledge of agent valuations. In dynamic V2V environments with private preferences and unpredictable departures, accurate valuation models are impractical, and the optimization must be re-solved from scratch as conditions change.

Multi-agent reinforcement learning (MARL) addresses model-based limitations by enabling agents to learn adaptive strategies through interaction \cite{qiu2023rl_review,park2022multi}. The centralized training with decentralized execution (CTDE) paradigm \cite{lowe2017maddpg} is well-suited: critics access joint information during training to handle non-stationarity, while actors use only local observations during deployment. MARL has been applied to EV charging recommendation \cite{zhang2021master}, auction-based EV charging coordination \cite{zou2022intelligent}, EV virtual power plant operation \cite{huang2026safe}, and microgrid trading with Shapley-value credit assignment \cite{harrold2022marl_trading}, where it outperformed limitations of single-agent methods. For V2V specifically, MARL with CTDE was applied for decentralized energy exchange in \cite{fan2023marl}, demonstrating that multi-agent approaches can coordinate V2V markets. However, the framework lacked game-theoretic incentive structures: agents optimized individual utility without fairness guarantees, and the market clearing had no formal Pareto efficiency or individual rationality properties, resulting in potentially arbitrary and unpredictable equilibria that may discourage participation.

The preceding review reveals a fundamental complementarity: NBS guarantees fairness but requires complete information unavailable in dynamic V2V settings, while MARL provides adaptive learning but converges to arbitrary equilibria without fairness assurance. We address this gap by using NBS simultaneously as market clearing rule and reward shaping signal. At the upper level, NBS-based bilateral pricing ensures Pareto-efficient clearing; at the lower level, NBS-derived price proximity rewards guide agents toward bargaining-optimal strategies during training. Unlike generic reward shaping, this mechanism inherits Pareto efficiency and individual rationality by construction, resolving the absence of fairness properties in prior V2V policies \cite{fan2023marl,harrold2022marl_trading}.

This paper proposes Nash-MADDPG for incentive-aligned V2V energy trading, formulated as bi-level optimization. We define \textit{incentive alignment} as: (i) \textit{individual rationality}, where no agent receives negative utility, (ii) \textit{Pareto efficiency}, where no alternative allocation improves one agent without harming another. Both properties are guaranteed by NBS \cite{nash1950bargaining,wang2016incentivizing}. 
The price proximity reward in Section~\ref{sec:evalution} penalizes deviations from Nash Bargaining prices during training, structurally biasing agents toward submitting prices closer to their value-consistent prices at convergence.

The main contributions are twofold: (i) We develop a V2V coordination framework for charging stations with continuous EV turnover, where the core engineering challenge is that vehicles arrive and depart unpredictably with private valuations, making static optimization infeasible. We address this by combining Nash Bargaining market clearing with adaptive MARL that remains effective under continuous agent turnover across varying EV population scales.
(ii) We integrate NBS with MADDPG through dual mechanisms that address distinct failure modes that arise in V2V markets with self-interested EVs: bilateral pricing prevents arbitrary equilibrium when agents hold private valuations under continuous turnover, while price proximity rewards resolve the absence of fairness guarantees when agents learn independently with coordination, together empirically improving social welfare, fairness, and training stability over auction-based baselines.

The remainder of this paper is organized as follows. Section~\ref{sec:problem} presents the V2V problem formulation. Section~\ref{sec:algorithm} details Nash-MADDPG. Section~\ref{sec:evalution} presents evaluation. Section~\ref{sec:conclusion} concludes this paper.

\section{System Model and Problem Formulation}\label{sec:problem}

\subsection{V2V Trading System and EV Model}

\begin{figure}[t]
    \centering
    \includegraphics[width=\columnwidth]{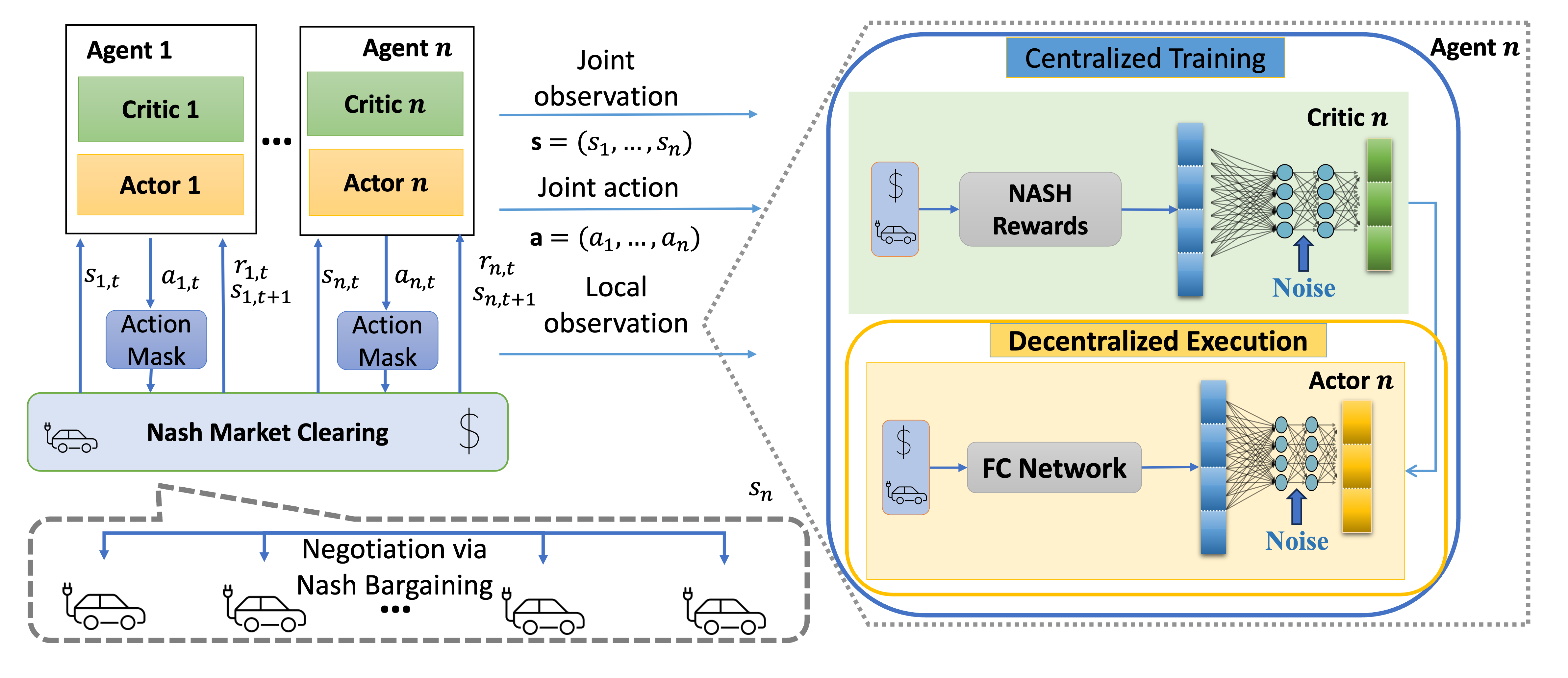}
    \caption{Nash-MADDPG system architecture showing a bi-level structure. Left: Multi-agent system with Nash Bargaining market clearing. Right: Centralized training and decentralized execution (CTDE) framework.}
    \label{fig:architecture}
\end{figure}

We consider a V2V energy trading system in a parking facility with dynamic EV population $\mathcal{V} = \{1, \ldots, N(t)\}$ operating over discrete timesteps $t \in \{0, 1, \ldots, T_{\text{max}}\}$ with 30-minute intervals (i.e., $\Delta t=0.5$). The system is structured as bi-level optimization: at the lower level, EV agents learn bidding policies $\{\pi_1, \ldots, \pi_N\}$ mapping local observations to price-quantity offers; at the upper level, market mechanism $\mathcal{M}$ clears trades to maximize Pareto efficiency under continuous agent turnover. Fig.~\ref{fig:architecture} illustrates the architecture.

Each EV $i$ maintains battery energy level $B_i(t)$ and specifies a target charge level $B^{\text{need}}_i$ upon arrival. Based on their energy state, EVs self-categorize as: \textit{buyers} if $B_i(t) < B^{\text{need}}_i$, \textit{sellers} if $B_i(t) > B^{\text{need}}_i + \Delta B^{\text{buffer}}$ (where $\Delta B^{\text{buffer}}$ is a safety margin), or \textit{neutral} otherwise. As time until departure decreases, urgency increases according to:
\begin{equation}
u_i(t) = u^{\text{base}}_i + \lambda_{\text{time}} \cdot \left(1 - \frac{T^{\text{remaining}}_i(t)}{T^{\text{duration}}_i(t)}\right),
\label{eq:urgency}
\end{equation}
where $u^{\text{base}}_i$ is baseline urgency, $\lambda_{\text{time}}$ scales urgency over time, $T^{\text{remaining}}_i(t)$ is time until departure, and $T^{\text{duration}}_i(t)$ is total parking duration.

Each buyer $i$'s willingness to pay is modeled as:
\begin{equation}
v_{b,i}(t) = v^{\text{base}}_b + \beta_{\text{urgency}} \cdot u_i(t) + \gamma_{\text{opportunity}} \cdot \frac{B^{\text{deficit}}_i(t)}{B^{\max}_i},
\label{eq:buyer_valuation}
\end{equation}
where $v^{\text{base}}_b$ is a baseline valuation, $\beta_{\text{urgency}}$ scales urgency impact on willingness to pay, $B^{\text{deficit}}_i(t) = B^{\text{need}}_i - B_i(t)$ is energy shortage, and $\gamma_{\text{opportunity}}$ weights the penalty for insufficient charge relative to battery capacity.

Each seller $j$'s minimum acceptable price is:
\begin{equation}
c_{s,j}(t) = c^{\text{base}}_s - v^{\text{battery}} \cdot \frac{B^{\text{avail}}_j(t)}{B^{\max}_j} + c^{\text{grid}} \cdot \rho_j + \delta^{\text{degrad}},
\label{eq:seller_cost}
\end{equation}
where $c^{\text{base}}_s$ is a baseline cost, $v^{\text{battery}}$ is the intrinsic value of stored energy, $B^{\text{avail}}_j(t) = B_j(t) - B^{\text{need}}_j - \Delta B^{\text{buffer}}$ is surplus energy, $c^{\text{grid}} \cdot \rho_j$ captures grid arbitrage opportunity with $\rho_j \sim \mathcal{U}(0, 1)$, and $\delta^{\text{degrad}}$ is per-kWh battery degradation cost. Feasibility is evaluated using submitted bids and asks before quantity allocation.

\subsection{Nash Bargaining Market Clearing}\label{sec:nash_clearing}

At each timestep, the coordinator partitions active agents into buyers $\mathcal{B}(t)$ and sellers $\mathcal{S}(t)$. True valuations are \textit{private}; the coordinator observes only submitted bids $b_i$ and asks $a_j$, and applies Nash Bargaining using these as effective valuations. For each feasible pair ($b_i \geq a_j$), the disagreement point is no-trade (zero utility for both), and the bilateral Nash Bargaining price maximizes the product of surpluses over disagreement \cite{nash1950bargaining}:
\begin{equation}
p^*_{ij} = \arg\max_p \; (b_i - p)(p - a_j),
\label{eq:nash_price}
\end{equation}
which, by first-order optimality with symmetric bargaining power, yields the closed-form midpoint:
\begin{equation}
p^*_{ij} = \tfrac{1}{2}(b_i + a_j).
\label{eq:nash_price_closed}
\end{equation}
Truth bidding is not a dominant strategy: by the Myerson-Satterwaite impossibility theorem~\cite{myerson1983efficient}, no bilateral-trade mechanism under private valuations can be simultaneously ex-post efficient, individually rational, budget-balanced, and incentive-compatible. Since our clearing is strongly budget-balanced (buyer payment equals seller revenue), under truthful bidding ($b_i = v_{b,i}$, $a_j = c_{s,j}$) the Pareto efficiency and individual rationality guarantees hold directly for true valuations, while strict incentive compatibility guarantees is necessarily relinquished. The price proximity in Eq.~\eqref{eq:price_proximity} drives learned policies toward value-consistent bidding so these guarantees hold approximately, with the residual gap quantified by the stable price band and individual-rationality results in Section~\ref{sec:evalution}.

Trade quantities are then determined by maximizing the log-transformed Nash social welfare over all feasible pairs. Let $\mathbf{X}(t) = \{x_{ij}(t)\}$, where $x_{ij}(t) \geq 0$ is energy transferred from seller $j$ to buyer $i$:
\begin{equation}
\begin{split}
& \mathbf{X}^*(t) \\
& = \arg\max_{\mathbf{X}(t)} \left[ \sum_{i \in \mathcal{B}} \log(U^{\text{total}}_{b,i}(t) + \epsilon) + \sum_{j \in \mathcal{S}} \log(U^{\text{total}}_{s,j}(t) + \epsilon) \right],
\label{eq:nash_product}
\end{split}
\end{equation}
where $\epsilon > 0$ prevents $\log(0)$, and agent utilities are:
\begin{align}
U^{\text{total}}_{b,i}(t) &= \sum_{j \in \mathcal{S}} (b_i(t) - p^*_{ij}(t)) \cdot x_{ij}(t), \\
U^{\text{total}}_{s,j}(t) &= \sum_{i \in \mathcal{B}} (p^*_{ij}(t) - a_j(t)) \cdot x_{ij}(t).
\label{eq:utilities}
\end{align}
For notational convenience, let $U_i(t)$ denote the total utility of agent $i$, 
i.e., $U_{b,i}^{\text{total}}(t)$ if $i \in \mathcal{B}$ and $U_{s,i}^{\text{total}}(t)$ if $i \in \mathcal{S}$.

The optimization is subject to: (i) $\sum_{i \in \mathcal{B}} x_{ij}(t) \leq \min(B^{\text{avail}}_j(t),\, 
q^{\text{offer}}_j(t),\, P^{\max}_j \cdot \Delta t)$ (seller capacity),
(ii) $\sum_{j \in \mathcal{S}} x_{ij}(t) \leq \min(B^{\text{deficit}}_i(t),\, 
q^{\text{offer}}_i(t),\, P^{\max}_i \cdot \Delta t)$ (buyer capacity), and (iii) $U^{\text{total}}_{b,i}(t) \geq 0$, $U^{\text{total}}_{s,j}(t) \geq 0$ (individual rationality). The log transformation converts the non-convex Nash product into a concave objective, forming a convex program solved via SLSQP. A single EV may trade with multiple counterparties per timestep; battery states update as $B_i(t+1) = B_i(t) + \eta \sum_{j} x^*_{ij}(t)$ for buyers and $B_j(t+1) = B_j(t) - \sum_{i} x^*_{ij}(t)$ for sellers, where $\eta \in (0,1]$ is charging efficiency.

\subsection{Multi-Agent MDP Formulation}\label{sec:mdp}

We formulate V2V coordination as a decentralized Multi-Agent MDP. Each agent $i$ observes local state $s_i(t) = \{B_i(t), B^{\text{deficit}}_i(t), B^{\text{avail}}_i(t), u_i(t), T^{\text{remaining}}_i(t), p^{\text{last}}_i(t), \text{role}_i(t)\}$ and outputs action $a_i(t) = [p_i(t), q_i(t)]$. Role-specific action masking enforces feasibility: neutral agents are excluded ($q_i = 0$), buyers are bounded by $q_i \in [0, B^{\text{deficit}}_i]$, and sellers by $q_i \in [0, B^{\text{avail}}_i]$. All agents share a single actor-critic network with role encoded as a one-hot feature, enabling generalization across unseen population sizes without retraining.

The information structure is asymmetric: true valuations $v_{b,i}$ and costs $c_{s,j}$ are \textit{private} and never revealed to the coordinator, which observes only submitted bids $b_i = p_i(t)$, asks $a_j = p_j(t)$, and quantities $q_i(t)$. Market clearing in Section~\ref{sec:nash_clearing} therefore uses submitted bids and asks as substitutes for true valuations. The price proximity reward in Eq.~\eqref{eq:price_proximity} encourages submitted prices toward Nash Bargaining prices during training, empirically reducing the gap between submitted and true valuations; to the extent that $b_i \approx v_{b,i}$ and $a_j \approx c_{s,j}$ at convergence, the individual rationality guarantees of NBS hold approximately from each agent's perspective.

Each agent maximizes expected cumulative discounted reward over its parking duration:
\begin{equation}
\max_{\pi_i} \mathbb{E}_{\pi_i} \left[ \sum_{t=t^{\text{entry}}_i}^{T^{\text{dep}}_i} \gamma^{t-t^{\text{entry}}_i} r_i(t) \right],
\label{eq:agent_objective}
\end{equation}
where $\gamma \in (0,1)$ is the discount factor. The reward $r_i(t)$ must be designed to align individual agent learning with system-level objectives.

The V2V system pursues three objectives: (i) Social Welfare $\text{(SW)} = \sum_{i \in \mathcal{B}} U^{\text{total}}_{b,i} + \sum_{j \in \mathcal{S}} U^{\text{total}}_{s,j}$, (ii) Jain's Fairness Index $FI = [\sum^K_{i=1} U_i]^2/(K \cdot \sum^K_{i=1} U^2_i)$ computed over $K$ active traders (buyers and sellers who submitted nonzero offers, excluding neutrals; zero-utility traders included) \cite{jain1984quantitative}, and (iii) Match Rate $P_{\text{match}} = |\{i : \sum_j x^*_{ij} > 0 \text{ or } \sum_k x^*_{ki} > 0\}| / (|\mathcal{B}| + |\mathcal{S}|)$, measuring the fraction of active participants successfully matched. Unlike pair-level agreement probability (which scales quadratically), $P_{\text{match}}$ is participant-level and comparable across population sizes. The challenge is designing $r_i(t)$ such that self-interested agents simultaneously achieve high performance on all three metrics.

\section{Nash-MADDPG Algorithm}\label{sec:algorithm}

V2V coordination faces three challenges: (i) non-stationary learning where strategies depend on simultaneously adapting opponents, (ii) pure learning converging to arbitrary outcomes without fairness properties, and (iii) continuous vehicle turnover requiring generalizable policies.

We propose Nash-MADDPG, integrating NBS into MADDPG through dual mechanisms. Nash Bargaining provides Pareto-efficient market clearing at the upper level; Nash-derived price proximity rewards guide learning at the lower level. The CTDE framework \cite{foerster2017stabilising} uses joint information during training and only local observations during deployment.

\subsection{Nash-MADDPG Architecture}

Each EV $i$ maintains actor network $\mu_i(s_i; \theta_i)$ and critic network $Q_i(s, a_1, \ldots, a_N; \phi_i)$, both implemented as multi-layer perceptrons with ReLU activations. The actor outputs tanh-scaled price and sigmoid-scaled quantity; the critic concatenates all agents' states and actions. Ornstein-Uhlenbeck noise provides exploration; states are normalized to $[0, 1]$.
\begin{equation}
Q_i(s, a_1, \ldots, a_N; \phi_i) = \mathbb{E} \left[ \sum_{t'=t}^{T_{\max}} \gamma^{t'-t} r_i(t') \right].
\label{eq:critic_value}
\end{equation}

The critic learns by minimizing prediction errors:
\begin{equation}
\mathcal{L}(\phi_i) = \mathbb{E}_{(s,a,r,s') \sim \mathcal{D}} \left[ (Q_i(s, a_1, \ldots, a_N; \phi_i) - y_i)^2 \right],
\label{eq:critic_loss}
\end{equation}
where $y_i = r_i + \gamma Q'_i(s', \mu'_1(s'_1), \ldots, \mu'_N(s'_N); \phi'_i)$ uses target networks $\mu'_i, Q'_i$ to stabilize learning \cite{lowe2017maddpg}. The actor updates by following the critic's gradient to maximize expected utility:
\begin{equation}
\nabla_{\theta_i} J(\theta_i) = \mathbb{E}_{s \sim \mathcal{D}} \left[ \nabla_{\theta_i} \mu_i(s_i; \theta_i) \nabla_{a_i} Q_i(s, a_1, \ldots, a_N; \phi_i) \right].
\label{eq:actor_gradient}
\end{equation}
Target networks update gradually with soft update parameter $\tau$ to maintain stability as vehicles arrive, depart, and update bidding strategies.

\subsection{Nash-Derived Reward Structure}

The reward aligns individual incentives with collective objectives (Section~\ref{sec:mdp}) through three components.

(1) \textit{Base utility}. The base utility $r^{\text{base}}_i(t) = U_i(p^*, q^*)$ captures immediate utility from executed trades at Nash-cleared prices $p^*$ and quantities $q^*$, maintaining individual rationality.

(2) \textit{Credit assignment}. Exact Shapley values are intractable, i.e., $O(2^N)$, so we use counterfactual baselines inspired by COMA \cite{foerster2018coma}: $\delta_i = r^{\text{collective}}(\mathbf{a}) - r^{\text{collective}}(\mathbf{a}_{-i}, \bar{a}_i)$, where $\bar{a}_i$ is the no-trade default, $\mathbf{a}_{-i}$ holds all other agents' actions fixed, and $r^{\text{collective}}(\mathbf{a}) = \lambda_w \cdot \text{SW}(\mathbf{a}) + \lambda_1 \cdot \text{FI}(\mathbf{a}) + \lambda_2 \cdot P_{\text{match}}(\mathbf{a})$ is a weighted sum of the three system objectives defined in Section~\ref{sec:mdp}, in which $\lambda_w$, $\lambda_1$, and $\lambda_2$ are weights. This $O(N)$ approximation maintains sensitivity to marginal contributions.

(3) \textit{Price proximity} (training only). The key innovation integrates bargaining-optimal pricing into learning during centralized training:
\begin{equation}
r^{\text{price}}_i(t) = \begin{cases} -\kappa_{\text{price}} \cdot |p_i(t) - p^{\text{Nash}}_{i\hat{j}}|^2 & \text{if } (i, \hat{j}) \in \mathbf{X}^* \\ 0 & \text{if unmatched}, \end{cases}
\label{eq:price_proximity}
\end{equation}
where $\hat{j} = \arg\max_{j: x^*_{ij} > 0} x^*_{ij}$ is the matched counterparty from $\mathbf{X}^*$, $\kappa_{\text{price}} > 0$ is the proximity weight, and $p^{\text{Nash}}_{i\hat{j}} = \tfrac{1}{2}(b_i + a_{\hat{j}})$ is computed from submitted bids/asks (not private valuations). This reward is computed by the coordinator during training only; at execution, actors use only local observations. The quadratic penalty biases policies toward bargaining-optimal prices, with Nash guidance propagating through critic training to actor updates.

Combining these three components, the complete reward for agent $i$ is:
\begin{align}
r_i(t) &= \alpha_r \cdot r^{\text{base}}_i(t) + \delta_i  + w_{\text{price}} \cdot r^{\text{price}}_i(t) \nonumber \\
&\quad - \mu_i \cdot \max\{0, -U_i(t)\},
\label{eq:integrated_reward}
\end{align}
where $\alpha_r$, $w_{\text{price}}$, and $\mu_i$ are tunable weights balancing individual utility, cooperation (scaled by $1/N$ for population-invariance), and Nash Bargaining guidance. The final penalty term discourages trades yielding negative utility, reinforcing individual rationality.

\section{Performance Evaluation}\label{sec:evalution}

\subsection{Experimental Settings}

We simulate Tesla Model 3 EVs (75 kWh capacity) with Australian electricity pricing (grid: \$0.28/kWh, V2V range: \$0.05-\$0.50/kWh). Training uses 2,000 episodes with $\gamma = 0.95$, actor learning rate $1 \times 10^{-4}$, critic learning rate $1 \times 10^{-3}$, replay buffer $10^5$, batch size 256, and $\tau = 0.01$. All agents share a single actor and critic network with role (buyer/seller) encoded as a one-hot input feature, enabling generalization across population sizes. Training uses $N = 20$ agents across three seeds; testing evaluates on $N \in \{6, 10, 15, 20, 30, 50, 75, 100\}$ without retraining. Testing includes one-day (16 timesteps) and 30-day (480 timesteps) evaluations, where the 30-day horizon encompasses multiple complete vehicle population turnovers, stress-testing policy generalization under continuous agent replacement.

The simulation environment is generated as follows. EV arrivals follow a Poisson process with rate $\lambda = N/8$ per timestep (calibrated to maintain target population $N$). Parking duration is drawn from $T^{\text{duration}}_i \sim \mathcal{U}(4, 12)$ timesteps (2-6 hours). Initial state-of-charge is $B_i(0)/B^{\max}_i \sim \mathcal{U}(0.2, 0.9)$, and target charge $B^{\text{need}}_i/B^{\max}_i \sim \mathcal{U}(0.5, 0.95)$. These parameters yield approximately 40-60\% buyers, 20-35\% sellers, and 15-30\% neutrals per timestep.

\subsection{Baseline Methods}

We compare Nash-MADDPG against three baselines. All methods use the same MADDPG agents outputting price-quantity actions $[p_i(t), q_i(t)]$; they differ only in the market clearing mechanism and reward structure.

(1) \textbf{Learning Only.} Standard MADDPG \cite{lowe2017maddpg} without Nash Bargaining clearing or fairness-shaped rewards; agents learn from immediate trade utility only.

(2) \textbf{Greedy Average.} The clearing price is set as the midpoint $p = \frac{1}{2}(b_i + a_j)$, identical to Nash Bargaining pricing, but quantities are allocated greedily by iterating over feasible pairs in random order without joint optimization. This isolates the contribution of Nash social welfare maximization in Eq.~\eqref{eq:nash_product} from the pricing rule, since both methods share the same bilateral price.

(3) \textbf{Double Auction.} buyer bids and seller asks are sorted, trades execute where buyer bid exceeds seller ask, and the clearing price is set at the bid-ask midpoint of the marginal pair.

\subsection{Results and Discussion}

\subsubsection{Learning Convergence}
\begin{figure}[t]
    \centering
    \includegraphics[width=\columnwidth]{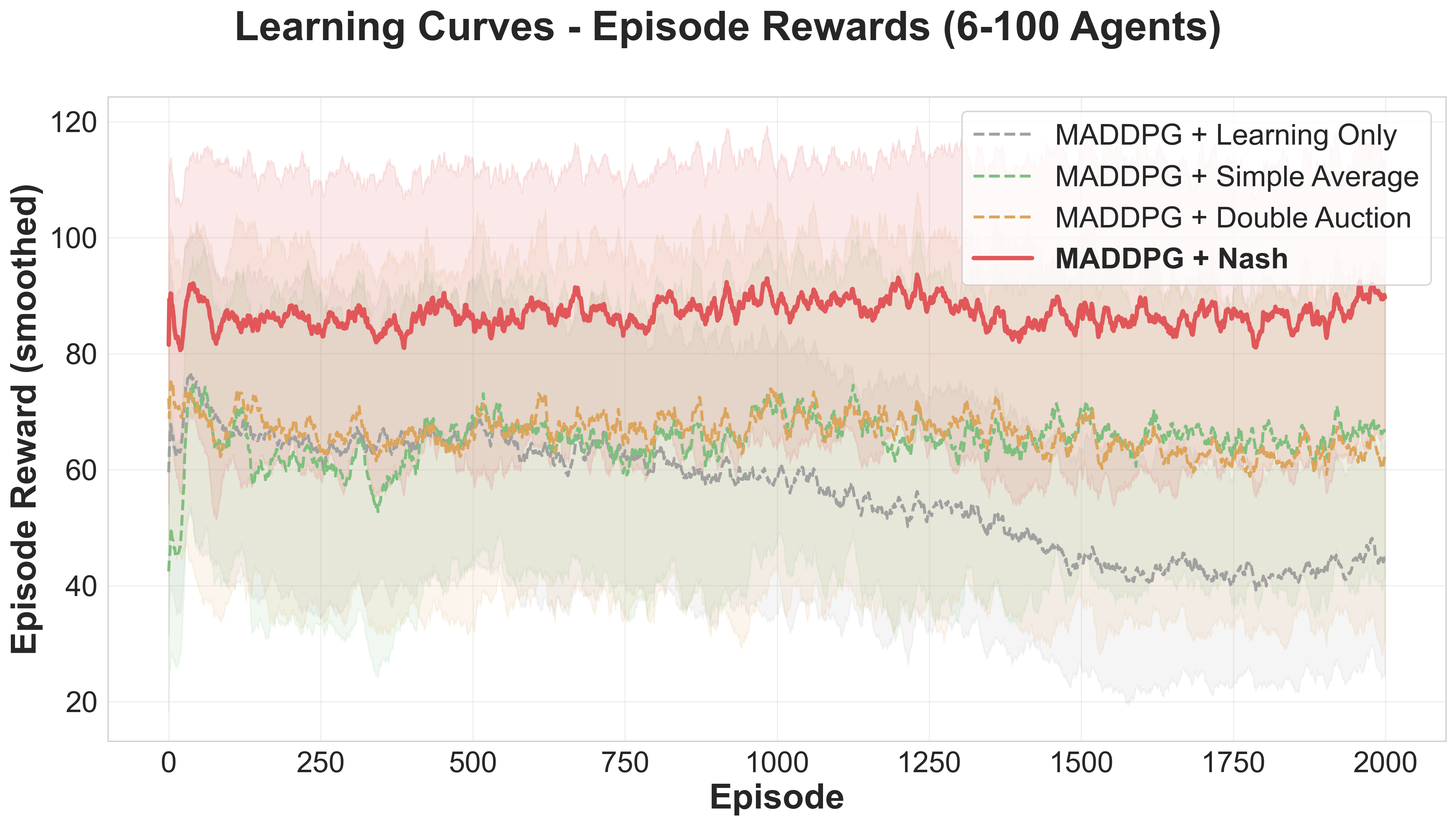}
    \caption{Learning curves aggregated across 6-100 agents over 2000 episodes. Nash-MADDPG maintains stable convergence near 90 throughout training.}
    \label{fig:learning_curves}
\end{figure}

Fig. \ref{fig:learning_curves} shows learning curves across 6-100 agents. Nash-MADDPG converges stably near 90 with minimal variance, while Learning Only collapses, and baselines plateau at 55-60 with higher variance.

\subsubsection{Time-Series Performance (30-Day Operation)}
Table \ref{tab:timeseries_results} presents comprehensive metrics from 30-day continuous operation aggregated across 6-100 agents. Nash-MADDPG outperforms other methods in terms of all metrics simultaneously.

\begin{table}[t]
    \centering
    \caption{Performance Metrics: 30 Days, 6-100 Agents (Median [IQR])}
    \label{tab:timeseries_results}
    \scriptsize
    \begin{tabularx}{\linewidth}{lXXXXX}
        \hline
        \textbf{Method} &
        \textbf{SW} &
        \textbf{Vol.} &
        \textbf{Gini} &
        \textbf{Jain's} &
        \textbf{$P_{\text{match}}$} \\
        \hline
        Learn. Only & 726 \scriptsize{[418]} & 4,504 \scriptsize{[2,581]} & .822 \scriptsize{[.247]} & .185 \scriptsize{[.135]} & .71 \scriptsize{[.18]} \\
        Greedy. Avg. & 607 \scriptsize{[289]} & 3,542 \scriptsize{[1,660]} & .797 \scriptsize{[.240]} & .216 \scriptsize{[.158]} & .65 \scriptsize{[.16]} \\
        Dbl. Auct. & 927 \scriptsize{[423]} & 6,037 \scriptsize{[2,832]} & .774 \scriptsize{[.233]} & .237 \scriptsize{[.173]} & .78 \scriptsize{[.15]} \\
        \textbf{Nash} & \textbf{1,497} \scriptsize{[264]} & \textbf{9,830} \scriptsize{[2,076]} & \textbf{.680} \scriptsize{[.126]} & \textbf{.332} \scriptsize{[.131]} & \textbf{.92} \scriptsize{[.07]} \\
        \hline
        \multicolumn{6}{l}{\scriptsize Median [IQR] over 3 seeds $\times$ 8 population sizes. SW (\$AUD),} \\
        \multicolumn{6}{l}{\scriptsize Vol (kWh). $P_{\text{match}}$: fraction of active participants matched.} \\
    \end{tabularx}
\end{table}

\begin{itemize}
    \item \textit{Economic efficiency}. In our experiments, Nash achieves 61.6\% improvement in social welfare (1,497 vs. 927 AUD) and 62.9\% in volume (9,830 vs. 6,037 kWh) over Double Auction, attributable to Nash social welfare allocation \cite{nash1950bargaining}.
    \item \textit{Fairness}. Jain's index is computed over active traders only (buyers and sellers with nonzero offers), including those matched with zero utility; neutral EVs are excluded. Nash achieves Gini 0.680 (vs. 0.774, lower is better) and Jain's 0.332 (vs. 0.237, higher is better), a 40.1\% improvement. 
    \item \textit{Trade completion}. Nash-MADDPG achieves $P_{\text{match}} = 0.92$ (vs. 0.78 for Double Auction), meaning 92\% of active participants are successfully matched, reflecting that individually rational pricing encourages broad participation.
\end{itemize}

\begin{figure}[t]
    \centering
    \includegraphics[width=0.94\columnwidth]{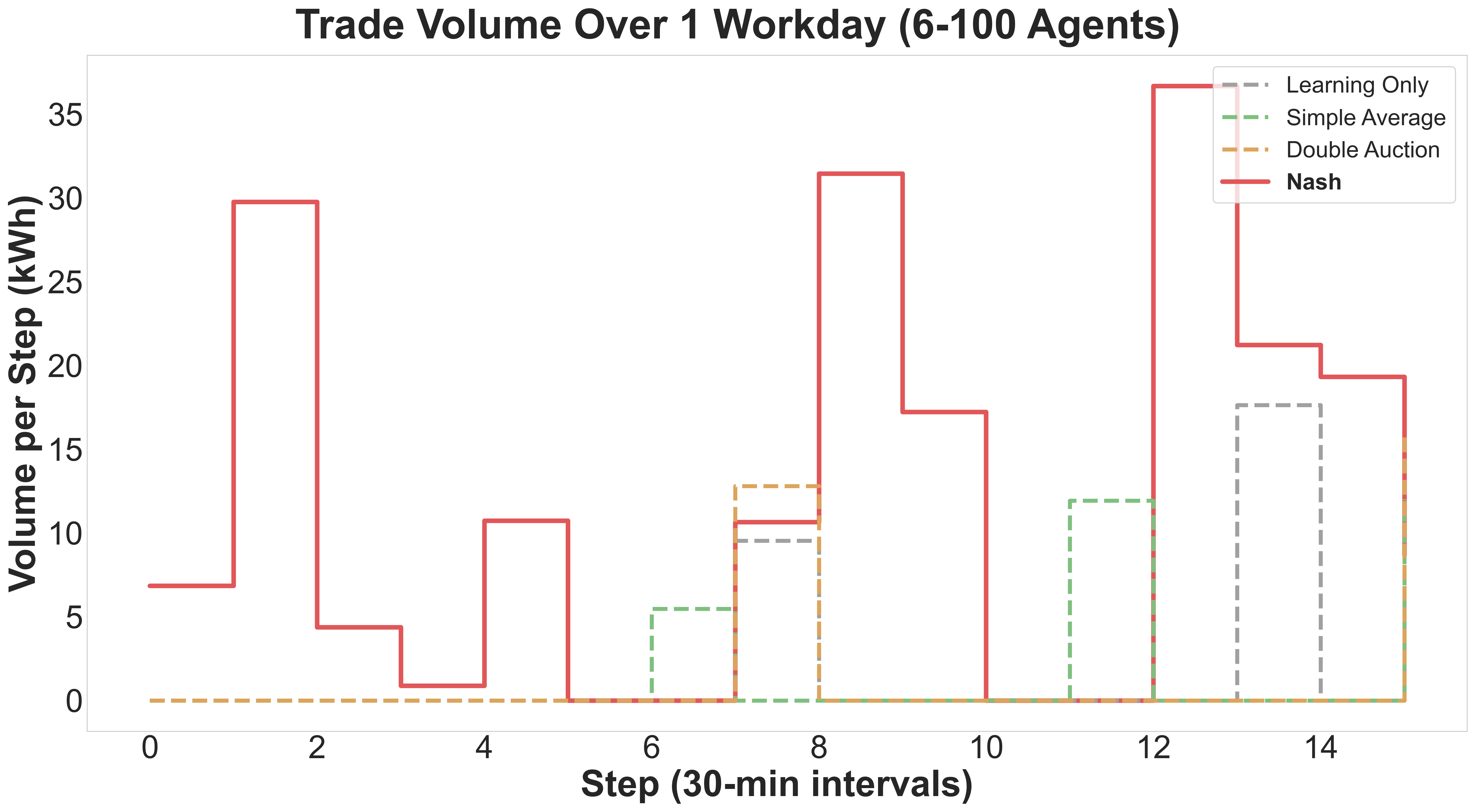}
    \caption{Trade volume over 8 hours of one workday. Nash-MADDPG adapts to intra-day demand dynamics with peak volumes during midday.}
    \label{fig:volume_evolution}
\end{figure}

\subsubsection{Intra-Day Market Dynamics}

Fig. \ref{fig:volume_evolution} shows Nash-MADDPG's three-phase workday pattern: morning ramp-up (steps 0-2, 30 kWh/step), midday peaks (steps 2-3, 8-12, up to 37 kWh/step), and afternoon decline (steps 12-16). Baselines show dampened responses (Learning Only peaks at 18 kWh, Double Auction at 13 kWh). 

\subsubsection{Long-Term Stability and Deployment Readiness}

The 480-timestep evaluation spans multiple complete vehicle population turnovers, demonstrating robust policy generalization under continuous turnover \cite{wang2023fast}. Nash-MADDPG achieves highest social welfare with fewer total trades (4,000 vs. 5,000 for Double Auction), indicating superior per-trade efficiency (0.75 vs. 0.55 utility/trade). Clearing prices exhibit stable convergence: mean 0.208 AUD/kWh (std 0.031), median 0.213 AUD/kWh, range [0.12, 0.34], positioned 25\% below grid price (0.28 AUD/kWh). 

\subsubsection{Scalability Robustness Across Population Sizes}

The 6-100 agent evaluation reveals consistency across population scales in our experiments. Nash-MADDPG achieves markedly lower coefficient of variation (CV = std/mean) than all baselines: 17.6\% for social welfare versus 45.6\% for Double Auction and 57.3\% for Learning Only, representing 2.6 times lower variability. Volume CV is similarly lower (21.1\% vs. 46.9\%), and fairness metrics show 36-43\% lower standard deviation. This consistency stems from Nash Bargaining's structural properties, whereas baselines rely on emergent coordination that degrades as action spaces grow.

\subsubsection{Effect of Nash Bargaining Guidance on Training Stability}
Fig.~\ref{fig:learning_curves} shows Learning Only undergoes catastrophic collapse with expanding confidence bands (20-80) illustrating non-stationarity. Nash-MADDPG maintains a tight band (83-95) across all configurations. Nash-MADDPG also activates trading from the first timestep, while baselines remain inactive until steps 6-7, representing approximately 3 hours of foregone trading.

\subsubsection{Reward Component Ablation}
Table~\ref{tab:ablation} isolates the contribution of each reward component by removing one term at a time from Eq.~\eqref{eq:integrated_reward}. Removing price proximity ($r^{\text{price}}$) causes the largest degradation: 34\% social welfare loss and training instability (reward variance increases 4.1 times), confirming it as the primary driver of both performance and stability. Removing credit assignment ($\phi_i$) reduces fairness most severely (Jain's drops 28\%), while removing global welfare/fairness terms ($\lambda_w, \lambda_1, \lambda_2$) has moderate impact on SW (12\%) but minimal effect on training stability.

\begin{table}[t]
    \centering
    \caption{Reward Ablation ($N = 100$, 30-Day Median)}
    \label{tab:ablation}
    \scriptsize
    \begin{tabularx}{\linewidth}{lXXXX}
        \hline
        \textbf{Variant} & \textbf{SW} & \textbf{Jain's} & \textbf{$P_{\text{match}}$} & \textbf{Stable?} \\
        \hline
        Full Nash-MADDPG & \textbf{1497} & \textbf{.332} & \textbf{.92} & \ding{51} \\
        $-\, r^{\text{price}}_i$ & 988 & .258 & .81 & \ding{55} \\
        $-\, \phi_i$ (credit) & 1301 & .239 & .87 & \ding{51} \\
        $-\, \lambda_w, \lambda_1, \lambda_2$ (global) & 1318 & .298 & .85 & \ding{51} \\
        \hline
        \multicolumn{5}{l}{\scriptsize Stable: reward variance $<$2 times of full model through ep. 2000.} \\
    \end{tabularx}
\end{table}

\section{Conclusion and Future Work}\label{sec:conclusion}

This paper proposes Nash-MADDPG for incentive-aligned V2V energy trading, integrating NBS into MADDPG through dual mechanisms: Pareto-efficient bilateral market clearing and Nash-derived price proximity rewards for learning guidance. Key design choices include log-transformed Nash social welfare for convex allocation, counterfactual credit assignment, and training-only reward shaping preserving decentralized execution. In our experiments across 6-100 agents over 30 days, Nash-MADDPG achieves 61.6\% higher social welfare, 62.9\% greater volume, and 40.1\% better fairness over Double Auction, with 2.6 times lower performance variability. Ablation confirms price proximity as the primary driver: removing it causes 34\% social welfare loss and training instability. These results suggest that game-theoretic integration substantially improves V2V market reliability under the studied conditions.

Future work will explore privacy-preserving federated training to eliminate centralized critic access, hierarchical coordination across multiple charging stations, and field validation with real EV charging data.

\bibliographystyle{IEEEtran}
\bibliography{ref}

\end{document}